\documentclass[12pt, leqno]{amsart}

\usepackage{shuffle} %\[A \shuffle B\]

\usepackage{lscape}
\usepackage{graphicx}
\usepackage{amssymb,amsmath,amsthm,amscd}
\usepackage{mathrsfs}
\usepackage{enumerate}
\usepackage[usenames,dvipsnames]{color}
\usepackage[colorlinks=true, pdfstartview=FitV,
 linkcolor=blue,citecolor=blue,urlcolor=blue]{hyperref}
\usepackage{upgreek}
\usepackage[all]{xy}
\allowdisplaybreaks[4]
\usepackage{oldgerm}
\usepackage{xspace}
\usepackage{marginnote}
\usepackage[normalem]{ulem}  % for strikeout \sout

\newcommand{\nc}{\newcommand}
\numberwithin{equation}{section}

\newenvironment{red}{\relax\color{red}}{\relax}
\newenvironment{blue}{\relax\color{Dandelion}}{\hspace*{.5ex}\relax}

\newcommand{\beb}{\begin{blue}}
\newcommand{\eb}{\end{blue}}

\newcommand{\berm}[1]{\begin{red}{}\marginnote{\fbox{\scshape\lowercase{M}}}%
#1}  % Masaki

\newcommand{\berE}[1]{\begin{red}{}\marginnote{\fbox{\scshape\lowercase{E}}}%
#1}  % Euiyong

\newcommand{\berMH}[1]{\begin{red}{}\marginnote{\fbox{\scshape\lowercase{MH}}}%
#1}  % Myungho

\newcommand{\berS}[1]{\begin{red}{}\marginnote{\fbox{\scshape\lowercase{S}}}
#1}  % Sejin

\renewcommand{\le}{\leqslant}
\renewcommand{\ge}{\geqslant}

\theoremstyle{plain}
\newtheorem{lemma}{Lemma}[section]
\newtheorem{prop}[lemma]{Proposition}
\newtheorem{theorem}[lemma]{Theorem}

\newcommand{\Prop}{\begin{prop}}
\newcommand{\enprop}{\end{prop}}
\newcommand{\Lemma}{\begin{lemma}}
\newcommand{\enlemma}{\end{lemma}}
\newcommand{\Th}{\begin{theorem}}
\newcommand{\enth}{\end{theorem}}
\newtheorem{corollary}[lemma]{Corollary}
\newcommand{\Cor}{\begin{corollary}}
\newcommand{\encor}{\end{corollary}}

\newtheorem{definition}[lemma]{Definition}
\newtheorem*{conjecture}{Conjecture}
\newcommand{\Def}{\begin{definition}}
\newcommand{\edf}{\end{definition}}
\newtheorem{sublemma}[lemma]{Sublemma}
\newcommand{\Sublemma}{\begin{sublemma}}
\newcommand{\ensub}{\end{sublemma}}

\theoremstyle{definition}
\newtheorem{remark}[lemma]{Remark}

\newtheorem{Convention}[lemma]{Convention}
\newcommand{\Conv}{\begin{Convention}}
\newcommand{\enconv}{\end{Convention}}
\nc{\Conj}{\begin{conjecture}}
\nc{\enconj}{\end{conjecture}}
\nc{\Rem}{\begin{remark}}
\nc{\enrem}{\end{remark}}
\newcommand{\C}{{\mathbb C}}

\newcommand{\Z}{{\mathbb Z}}
\newcommand{\B}{{\mathbf{B}}}

\newcommand{\sfW}{{\mathsf W}}

\newcommand{\DC}{{\rm D}}

\newcommand{\D}{\mathscr{D}}

\newcommand{\seteq}{\mathbin{:=}}

\newcommand{\hd}{{\operatorname{hd}}}

\newcommand{\g}{{\mathfrak{g}}}

\newcommand{\M}{{\mathcal M}}

\newcommand{\eq}{\begin{eqnarray}}
\newcommand{\eneq}{\end{eqnarray}}

\newcommand{\hs}{\hspace*}
\newcommand{\ms}{\mspace}

\newcommand{\To}[1][{\hs{2ex}}]{\xrightarrow{\,#1\,}}

\newcommand{\eqn}{\begin{eqnarray*}}
\newcommand{\eneqn}{\end{eqnarray*}}
\newcommand{\on}{\operatorname}

\newcommand{\sotimes}{\mathop{\mbox{\normalsize$\bigotimes$}}\limits}

\newcommand{\ba}{\begin{array}}
\newcommand{\ea}{\end{array}}

\newcommand{\set}[2]{\left\{#1 \mid #2 \right\}}

\newcommand{\eqsub}{\begin{subequations}\begin{eqnarray}}
\newcommand{\eneqsub}{\end{eqnarray}\end{subequations}}

\newcommand{\ol}{\overline}

\nc{\la}{\lambda}
\nc{\lam}{\lambda}
\nc{\U}[1][\g]{U_q(#1)}
\nc{\te}{\tilde{e}}
\nc{\tei}{\tilde{e}_i}
\nc{\tf}{\tilde{f}}
\nc{\tfi}{\tilde{f}_i}
\nc{\tU}{\widetilde U_q(\g)}
\nc{\tE}{\tilde{E}}
\nc{\tF}{\widetilde{\F}}
\nc{\tK}{\widetilde{K}}

\nc{\tk}{\tilde{k}}
\nc{\tkone}{\tk_{\ol{1}}}
\nc{\teone}{\tilde{e}_{\ol{1}}}
\nc{\tfone}{\tilde{f}_{\ol{1}}}

\nc{\teibar}{\tilde{e}_{\ol{i}}} \nc{\tfibar}{\tilde{f}_{\ol{i}}}
\nc{\tki}{{\tk}_{\ol {i}}}

\nc{\BZ}{{\mathbb{Z}}}
\nc{\al}{\alpha}
\nc{\qs}{{q}}
\nc{\lan}{\langle}
\nc{\ran}{\rangle}
\nc{\re}{{\mathrm{re}}}
\nc{\wt}{\operatorname{wt}}
\nc{\ch}{\operatorname{ch}}
\nc{\Um}[1][\g]{U^-_q(#1)}
\nc{\Ue}{U^+_q(\g)}
\nc{\eps}{\varepsilon}
\nc{\vphi}{\varphi}
\nc{\sphi}{\varphi^*}
\nc{\seps}{\varepsilon^*}

\nc{\nn}{\nonumber}
\def\max{{\mathop{\mathrm{max}}}}
\nc{\vp}{\varpi}
\nc{\cls}{{\operatorname{cl}}}
\nc{\Wt}{{\operatorname{Wt}}}
\nc{\Us}{U'_q(\g)}
\nc{\La}{\Lambda}
\nc{\tLa}{\widetilde\Lambda}
\nc{\ro}{{\rm(}}
\nc{\rf}{{\rm)}}
\nc{\norm}{{\mathrm{norm}}}
\nc{\qbox}{\quad\mbox}
\nc{\braid}{{\mathfrak{B}}}
\nc{\Ad}{\operatorname{Ad}}
\nc{\Aut}{\operatorname{Aut}}
\nc{\dt}[1]{\tilde{\tilde #1}}
\nc{\Sn}{S^{{\mathrm{norm}}}}
\nc{\aff}{{\rm{aff}}}
\nc{\rk}{{\mathrm{rk}}}
\nc{\tP}{\widetilde{P}}
\nc{\tW}{\widetilde{W}}
\nc{\Dyn}{\mathrm{Dyn}}
\nc{\tD}{\widetilde{\Delta}}
\nc{\height}[1]{{\operatorname{ht}}(#1)}
\nc{\bl}{\bigl(}
\nc{\br}{\bigr)}
\nc{\Hecke}{\mathrm{H}}
\nc{\HA}{\Hecke^{\mathrm{A}}}
\nc{\HB}{\Hecke^{\mathrm{B}}}
\newcommand{\scbul}{{\,\raise1pt\hbox{$\scriptscriptstyle\bullet$}\,}}
\nc{\vac}{{\phi}}
\nc{\Bt}{\B_\theta(\g)}
\nc{\be}{\begin{enumerate}}
\nc{\ee}{\end{enumerate}}
\nc{\low}{{\mathrm{low}}}
\nc{\upper}{{\mathrm{up}}}
\nc{\Zodd}{\Z_{\mathrm{odd}}}
\nc{\Ft}[1][n]{\mathbb{P}\mathrm{ol}_{#1}}
\nc{\Ftf}[1][n]{\widetilde{\mathbb{P}\mathrm{ol}}_{#1}}
\nc{\KA}{\on{K}^{\mathrm{A}}}
\nc{\KB}{\on{K}^{\mathrm{B}}}
\nc{\Res}{\on{Res}}
\nc{\Fc}[1][{n,m}]{\mathbf{F}_{#1}}
\nc{\tphi}{\tilde{\varphi}}
\nc{\CO}{\mathscr{O}}
\nc{\inte}{\mathrm{int}}
\nc{\Oint}{\mathcal{O}^{\ge0}_{\inte}}
\nc{\vs}{\vspace*}
\nc{\tLt}{\widetilde{L}}
\nc{\tL}{\widetilde{\Lambda}}
\nc{\tu}{\tilde{u}}
\nc{\noi}{\noindent}
\nc{\heigh}{\mathfrak{t}}
\nc{\lowest}{\mathfrak{l}}
\nc{\rootl}{\mathsf{Q}}
\nc{\cl}{\colon}
\nc{\uqpg}{U'_q(\mathfrak g)}
\nc{\uq}{\uqpg}
\nc{\Oh}{\widehat{\mathcal{O}}}
\nc{\pn}{p_{\mathfrak{n}}}

\nc{\KLR}{KLR algebra}
\nc{\tckname}{envelope} %{envelope}
\nc{\tclname}{range} %{area}

\nc{\KLRs}{KLR algebras}
\nc{\cor}{\mathbf{k}}
\nc{\cora}{{\cor(A)}}
\nc{\haut}{\mathrm{ht}}
\nc{\tens}{\mathop\otimes}
\nc{\gmod}{\mbox{-$\mathrm{gmod}$}}
\nc{\gMod}{\mbox{-$\mathrm{gMod}$}}
\nc{\proj}{\mbox{-$\mathrm{proj}$}}
\nc{\gproj}{\mbox{-$\mathrm{gproj}$}}
\nc{\smod}{\mbox{-$\mathrm{mod}$}}
\nc{\Mod}{\mbox{-$\mathrm{Mod}$}}
\nc{\h}{\mathfrak h}
\nc{\Rnorm}{R^{\rm{norm}}}
\nc{\Runiv}{R^{\rm{univ}}}
\nc{\Rren}{R^{\rm{ren}}}

\nc{\Vhat}{\widehat{V}}
\nc{\F}{\mathcal{F}}

\def\T{{\mathcal T}}

\nc{\fd}[1][A]{\on{\mathrm{flat.dim}_{#1}}}
\nc{\bP}{{\mathbb{P}}}
\nc{\bPh}{\widehat{\mathbb{P}}}
\nc{\bK}[1][{n}]{\widehat{\mathbb{K}}_{#1}}
\nc{\bV}[1][{n}]{\widehat{V}^{\otimes{#1}}}
\nc{\bVK}[1][{n}]{\widehat{V}^{\otimes{#1}}_{\widehat{\mathbb{K}}}}
\nc{\hV}{\widehat{V}}
\nc{\opp}{\mathrm{opp}}
\nc{\col}{\colon}
\nc{\bnum}{\be[{\rm(i)}]}
\nc{\oep}{\epsilon}
\nc{\qtext}{\quad\text}
\nc{\qtextq}[1]{\quad\text{#1}\quad}
\nc{\longtwoheadrightarrow}[1][]{\xymatrix{\ar@{->>}[r]^-{{#1}}&}}
\nc{\epiTo}[1][]{\longtwoheadrightarrow[{#1}]}
\nc{\epito}{\twoheadrightarrow}
\nc{\monoTo}[1][]{\xymatrix{\ar@{>->}[r]^-{{#1}}&}}
\nc{\sym}{\mathfrak{S}}
\nc{\inp}[1]{{({#1})_{\mathrm{n}}}}
\nc{\rtl}{\rootl}
\nc{\wtd}{\widetilde}
\nc{\etens}{\boxtimes}
\nc{\ds}[1]{\mathrm{d}(#1)}
\nc{\rmat}[1]{{\mathbf{r}}_%
{\mspace{-2mu}\raisebox{-.6ex}{${\scriptstyle{#1}}$}}}
\nc{\rmats}[1]{{\mathbf{r}}_%
{\mspace{-2mu}\raisebox{-.6ex}{${\scriptscriptstyle{#1}}$}}}
\nc{\shc}{\mathcal{C}}
\nc{\shs}{\mathcal{S}}
\nc{\Fct}{{\on{Fct}}}
\nc{\tC}{\widetilde{\shc}}
\nc{\Zp}{\Z_{\ge0}}
\nc{\tPhi}{\widetilde{\Phi}}
\nc{\tT}{{\widetilde{\T}}}
\nc{\Ob}{\on{Ob}}
\nc{\bwr}{\mbox{\large$\wr$}}
\nc{\Img}{\on{Im}}
\nc{\Ab}{\mathcal{A}^{\mathrm{big}}}
\nc{\Sb}{\mathcal{S}^{\mathrm{big}}}
\nc{\As}{\mathcal{A}}
\nc{\Ss}{\mathcal{S}}
\nc{\ntens}{\widetilde{\otimes}}
\nc{\hR}{\widehat{R}}
\nc{\nconv}{\mathop{\mbox{\large $\odot$}}}
\nc{\snconv}{\mbox{\scriptsize$\odot$}}
\nc{\ts}{\tilde{s}}
\nc{\sho}{\mathcal{O}}
\nc{\bc}{\begin{cases}}
\nc{\ec}{\end{cases}}
\nc{\slnh}{{\widehat{\mathfrak{sl}}_N}}
\nc{\UA}{U_q'(\slnh)}
\nc{\KR}{R_K}
\nc{\cQ}{\mathcal{Q}}
\nc{\rmQ}{\mathrm{Q}}
\nc{\Irr}{\mathcal{I}rr}
\nc{\tQ}{\widetilde{\cQ}}
\nc{\bs}{\mathbf{s}}
\nc{\bL}{\mathbb{L}}
\nc{\tg}{\tilde{g}}

\nc{\conv}{\mathbin{\mbox{\large $\circ$}}}
\nc{\shconv}{\mathbin{\large\diamond}}
\nc{\sconv}{\mathbin{\large\Delta}}
\nc{\hconv}{\mathbin{\nabla}}
\nc{\Rm}{R^{\mathrm{ren}}}

\nc{\bQ}{\ol{Q}}

\nc{\de}{\on{\textfrak{d}\ms{1mu}}}

\nc{\xmono}{\ar@{>->}}
\nc{\xepi}{\ar@{->>}}
\nc{\db}[1]{\raisebox{-.5ex}[2ex][1.8ex]{$#1$}}
\nc{\wb}[1]{\mbox{$\rule[-1.1ex]{0ex}{2ex}#1$}}
\nc{\univ}{\mathrm{univ}}
\nc{\rM}{{}^*\mspace{-2mu}M}
\nc{\lM}{M^*}
\nc{\uqm}{\uq\smod}
\nc{\tR}{\widetilde{R}_{\gamma,\beta}}
\nc{\tx}{\tilde{x}}
\nc{\bi}{\mathbf{i}}
\nc{\ttau}{\widetilde{\tau}}

\nc{\tEnd}{\on{\widetilde{E}nd}}
\nc{\tHom}{\on{\widetilde{H}om}}

\nc{\K}{{K}}
\nc{\Kex}{{\K}_{\mathrm{ex}}}
\nc{\Kfr}{{\K}_{\mathrm{f\mspace{.01mu}r}}}
\nc{\coro}{\cor}
\nc{\tB}{\widetilde{B}}
\nc{\seed}{\mathscr{S}}
\nc{\up}{\mathrm{up}}
\nc{\bfa}{\mathbf{a}}

\newlength{\mylength}
\nc{\ov}[1]{\overline{#1}}
\nc{\Wlmj}[3]{\W_{#2,#3}^{(#1)}}
\nc{\Mkl}[2]{\M_\ttww(#1,#2)}
\nc{\mqs}{(-q^2)}
\nc{\Cquiver}{\upsigma}

\nc{\mut}[1]{{\mu}_{\mspace{-2mu}\raisebox{-.5ex}{${\scriptstyle{#1}}$}}}

\nc{\Kt}{\mathcal K_t}
\nc{\KT}{\mathbb{K}_t}
\nc{\yim}{y_{i,m}}
\nc{\yjm}{y_{j,m}}
\nc{\yjp}{y_{j,p}}
\nc{\yimp}{y_{i,m+1}}
\nc{\yjmp}{y_{j,m+1}}

\nc{\Refl}{\mathscr{S}}
\nc{\Reflinv}{{\Refl}^{-1}}

\nc{\catC}{\mathscr C}
\nc{\catA}{\mathcal A}
\nc{\scrA}{{\mathscr{A}}}
\nc{\frakC}{{\mathfrak{C}}}
\nc{\frakT}{{\mathfrak{T}}}
\nc{\shift}{{\mathrm T}}
\nc{\rE}{ \mathsf{E} }
\nc{\rW}{ \mathcal{W} }
\nc{\rES}{ \mathcal{E} }

\nc{\brd}{\sigma} %generator of the braid group
\nc{\into}{\xymatrix@C=3ex{{}\ar@{^{(}->}[r]&{}}}
\nc{\dual}{\D}
\nc{\cat}[1][{\g}]{\catC_{#1}^0}
\nc{\qt}[1]{[{#1}]_t}

\nc{\catCO}{{\catC_\g^0}}
\nc{\catCQ}{{\catC_{\qQ}}}
\nc{\Li}{{\La^\infty}}
\nc{\sigZ}{{\sigma_0(\g)}}
\nc{\sigQ}{{\sigma_\qQ(\g)}}

\nc{\ZZ}{{\mathbf{Z}}}
\nc{\sP}{{\mathsf{P}}}
\nc{\sV}{{\mathsf{V}}}
\nc{\rxw}{{\underline{w_0}}}
\nc{\boten}[1]{\overrightarrow{\bigotimes_{#1}}}
\nc{\cmA}{{\mathsf{A}}}
\nc{\cmC}{{\mathsf{C}}}
\nc{\ddD}{{\mathcal{D}}}
\nc{\qQ}{{\mathcal{Q}}}
\nc{\gf}{{\g_{\rm fin}}}
\nc{\If}{{I_{\rm fin}}}
\nc{\cmAf}{{\cmA_{\rm fin}}}
\nc{\weyl}{{\mathsf{W}}}
\nc{\weylf}{{\weyl_{\rm fin}}}
\nc{\Deg}{\mathrm{Deg}}
\nc{\KRc}{{K_{q=1}(R_\cmC\gmod)}}
\nc{\prD}{{\Delta^+_{0}}}
\nc{\n}{{\mathfrak{n}}}
\nc{\finn}{{\mathfrak{n}}}
\nc{\Rt}{L} %simple root modules in the duality datum
\nc{\Cp}{V} %cuspidal modules
\nc{\cuspS}{{\mathsf{S}}}
\nc{\st}[1]{\left\{{#1}\right\}}
\nc{\WS}{quantum affine Weyl-Schur duality\xspace}
\nc{\CWS}{Quantum affine Weyl-Schur duality}
\nc{\hw}{\widehat{w}}
\newcommand{\tc}{{\widetilde{c}}}
\newcommand{\tb}{{\widetilde{b}}}
\newcommand{\ta}{{\widetilde{a}}}
\newcommand{\LL}{\mathcal{L}}
\newcommand{\RR}{\mathcal{R}}
\newcommand{\hI}{{\widehat{I}_0}}
\nc{\cd}[1]{\left\vert\ms{1mu}{#1}\ms{1mu}\right\vert}
\nc{\qtq}[1][{and}]{\quad\text{#1}\quad}

\title{Categories over quantum affine algebras and monoidal categorification}

\author[M. Kashiwara]{Masaki Kashiwara}
\thanks{The research of M.\ Kashiwara
was supported by Grant-in-Aid for Scientific Research (B)
15H03608, Japan Society for the Promotion of Science.}
\address[M. Kashiwara]{
Kyoto University Institute for Advanced Study,
Research Institute for Mathematical Sciences, Kyoto University,
Kyoto 606-8502, Japan \& Korea Institute for Advanced Study, Seoul 02455, Korea }
\email{masaki@kurims.kyoto-u.ac.jp}

\author[M. Kim]{Myungho Kim}
\address[M. Kim]{Department of Mathematics, Kyung Hee University, Seoul 02447, Korea}
\email{mkim@khu.ac.kr}
\thanks{The research of M.\ Kim was supported by the National Research Foundation of
Korea(NRF) Grant funded by the Korea government(MSIP) (NRF-2017R1C1B2007824).}

\author[S.-j. Oh]{Se-jin Oh}
\thanks{ The research of S.-j.\ Oh was supported by the Ministry of Education of the Republic of Korea and the National Research Foundation of Korea (NRF-2019R1A2C4069647).}
\address[S.-j. Oh]{Department of Mathematics, Ewha Womans University, Seoul 03760, Korea}
\email{sejin092@gmail.com}

\author[E. Park]{Euiyong Park}
\address[E. Park]{Department of Mathematics, University of Seoul, Seoul 02504, Korea}
\email{epark@uos.ac.kr}

\keywords{ {Monoidal categorification}, {Quantum affine algebra},
 {Cluster algebra}, {Kirillov-Reshetikhin module}, {$T$-system}}
\subjclass[2010]{17B37, 81R50, 18D10 } %

\date{May 22, 2020}

\begin{document}
\maketitle

\begin{abstract}
Let $U_q'(\g)$ be a quantum affine algebra of untwisted affine $ADE$ type,
and $\catC_\g^0$ the Hernandez-Leclerc category
of finite-dimensional $U_q'(\g)$-modules.
For a suitable infinite sequence $\widehat{w}_0= \cdots   s_{i_{-1}}s_{i_0}s_{i_1} \cdots$ of simple reflections,
we introduce subcategories $\catC_\g^{[a,b]}$ of $\catC_\g^0$ for all $a \le b \in \Z \sqcup\{ \pm \infty \}$.
Associated with a certain chain $\mathfrak{C}$ of intervals in $[a,b]$,
we construct a real simple commuting family $M(\mathfrak{C})$
in $\catC_\g^{[a,b]}$, which consists of Kirillov-Reshetikhin modules.
The category $\catC_\g^{[a,b]}$ provides a monoidal categorification of the cluster algebra $K(\catC_\g^{[a,b]})$, whose set of initial cluster variables is $[M(\mathfrak{C})]$.
In particular, this result gives an affirmative answer to the monoidal categorification conjecture on $\catC_\g^-$ by Hernandez-Leclerc since it is $\catC_\g^{[-\infty,0]}$, and is also applicable to
$\catC_\g^0$ since it is $\catC_\g^{[-\infty,\infty]}$.
\end{abstract}

\section{Introduction}
Let $U_q'(\g)$ be a quantum affine algebra.
The category $\catC_\g$ of finite-dimensional
integrable modules over $U_q'(\g)$ has been intensively studied due to its rich structure. For instances, every object $M$ in $\catC_\g$ has its left $M^*$ and right dual ${}^*M$, and
the $q$-characters of Kirillov-Reshetikhin modules in $\catC_\g$ provide a solution of the $T$-system, a system of differential equations appearing in solvable lattice models (\cite{FR99,Her06,KNS,Nak01}).

On the other hand, the cluster algebras were introduced by Fomin and Zelivinsky in~\cite{FZ02} to investigate upper global bases and total positivity in an aspect of combinatorics.

Interestingly, it is proved in~\cite{HL10,HL11,HL16} that the Grothendieck rings $K(\shc)$ of  monoidal subcategories $\shc=\catC_N$ $(N \in \Z_{\ge1})$, $\catC_\cQ$, $\catC_\g^{-}$  of $\catC_\g$ have cluster algebra structures $\scrA$, and conjectured that
every cluster monomial corresponds to the isomorphism class of a \emph{real} simple module in $\shc$; that is, $\shc$ is expected to be a \emph{monoidal categorification} of $\scrA$.
The conjectures for $\catC_N$ $(N \in \Z_{\ge1})$
of untwisted affine $ADE$ types are
proved in  \cite{HL10,HL13,Nak11}  and \cite{Qin17}.
Also, the conjecture for the subcategory $\catC_\cQ \subset \catC_\g$, determined by a \emph{$\rmQ$-data} $\cQ=(Q,\phi_Q)$ (\cite{FO20,KKOP20B}), is proved in~\cite{KKKO18} via the quantum affine Weyl-Schur duality functor $\F_\cQ$
(\cite{KKK18B,KKKO16D,KO18,OT19})
from the category $\shc_{QH}$ of finite-dimensional graded modules over the symmetric quiver Hecke algebra to $\catC_\cQ$.
More precisely, the category $\shc_{QH}$
provides a monoidal categorification of the quantum cluster algebra
$A_q(\finn)$, the quantum unipotent coordinate algebra of finite simply-laced type (\cite{BZ05}).
Since $\F_\cQ$ is an exact monoidal functor preserving simplicity,
we can prove the conjecture for $\catC_\cQ$ in an indirect way. However,
this method could not be applicable to other $\shc$ directly.

Recently, in~\cite{KKOP19C}, the authors of the present paper (KKOP)
developed $\Z$-valued invariants $\La,\La^\infty,\tilde{\La},\de$ for pairs of modules in $\catC_\g$, which is extracted from distinguished $U_q'(\g)$-module
homomorphisms, called \emph{$R$-matrices}. Furthermore, KKOP provided a criterion for a monoidal subcategory $\catC \subset \catC_\g$ to become a monoidal categorification of a cluster algebra
by using those invariants.
This paper can be understood as a continuation of~\cite{KKOP19C}, since we will apply the above criterion to various subcategories $\catC$ of $\catC_\g$,
including $\catC_\g^0$, $\catC_\g^-$ and $\catC_N$.
We also give their initial monoidal seeds in a uniform manner.

\smallskip
Let $\g_0$ be a finite-dimensional simple Lie algebra of $ADE$ type with a Cartan matrix $\cmA=(a_{ij})_{i,j\in I_0}$, $\sfW$ the Weyl group
generated by simple reflections $s_i$ $(i \in I_0)$, $\g$ the untwisted affine Kac-Moody algebra associated with $\g_0$, and $\uqpg$ the quantum affine algebra associated with $\g$.
In~\cite{HL10},
Hernandez-Leclerc defined the full subcategory $\catC_\g^0$ of $\catC_\g$.
Since every simple modules in $\catC_\g$ is a tensor product of suitable parameter shifts of simple modules in $\catC_\g^0$, it is enough to consider subcategories of  $\catC_\g^0$.

By extending a reduced expression $s_{i_1}  s_{i_{2}}\cdots s_{i_\ell}$ of the longest element $w_0$ of the Weyl group $\sfW$,
we obtain an infinite sequence
\begin{align} \label{eq: hw_0 intro}
\widehat{w}_0= \cdots  s_{i_{-2}}  s_{i_{-1}}s_{i_0}s_{i_1}  s_{i_{2}} \cdots
\end{align}
 of simple reflections satisfying properties \eqref{it:i} and \eqref{it:ii} in Section~\ref{sec: subcategory},
and then we define fundamental modules
$V[k]^{\widehat{w}_0}$ ($k\in\Z$).
For each \emph{interval} $[a,b]=\{ k \in \Z \mid a \le k \le b \}$ with $a \le b \in \Z \sqcup \{ \pm \infty \}$,
we define the subcategory $\catC_\g^{[a,b]}$ of $\catC_\g^0$  which is the smallest full monoidal subcategory containing $V[k]^{\widehat{w}_0}$ for all $k \in [a,b]$.
Then $\catC_\g^0$ is nothing but $\catC^{[-\infty,+\infty]}_\g$ and
the subcategory $\catC_\g^-$ introduced by Hernandez-Leclerc (\cite{HL16})
can be identified with $\catC_\g^{[-\infty,0]}$ (Remark~\ref{rmk: knownC}).

We say that  an interval $[a,b]$ is an \emph{$i$-box}
if ${i_{a}}={i_{b}}$. For each $i$-box $[a,b]$, we define a simple module $M[a,b]$, which can be understood as a \emph{quantum affine analogue of the determinantial module} (see Remark~\ref{rmk: determinantial module}).
In Theorem~\ref{thm: M[a,b]}, we show that $M[a,b]$ is a Kirillov-Reshetikhin module and give a sufficient condition for
the simplicity of the tensor product $M[a,b] \tens M[a',b']$
for $i$-boxes $[a,b]$ and $[a',b']$.
Then we define the notion of an \emph{admissible chain $\frakC=\{  [a_k,b_k] \subset [a,b] \mid 1 \le k\le b-a+1, \;{i_{a_{k}}}={i_{b_{k}}} \}$ of $i$-boxes} for an interval $[a,b]$ satisfying certain properties
(Definition~\ref{def: admissible chain}).
For each admissible chain $\frakC$, the family of Kirillov-Reshetikhin modules $M(\frakC)=\{ M[a_k,b_k] \}_{1\le k\le b-a+1}$ in $\catC_\g^{[a,b]}$ forms
a commuting family of real simple modules (Theorem~\ref{thm: real commuting family}).

\smallskip

The next step is to show that
$K(\catC_\g^{[a,b]})$ has a cluster algebra structure, $\catC_\g^{[a,b]}$ provides
a monoidal categorification of $K(\catC_\g^{[a,b]})$, and
any admissible chain $\frakC$ gives a monoidal seed $M(\frakC)$ (Theorem~\ref{thm:main}).
Based on the criterion in \cite{KKOP19C}, we shall prove this by showing the assertion for a special chain $\frakC$ , and then by extending it to a general
$\frakC$. Namely, we proceed by proving
\begin{enumerate}[{\rm (i)}]
\item the existence of a $\Uplambda$-admissible monoidal seed $\seed$ of $K(\catC_\g^{[a,b]})$ whose set of
initial cluster variable modules is $M(\frakC)$ for some admissible chain $\frakC$,
\item the existence of sequences of mutations among the $M(\frakC)$'s
 \emph{only employing} $T$-systems,
\end{enumerate}
which implies that any admissible chain $\frakC$ gives a $\Uplambda$-admissible monoidal seed for
all $\catC_\g^{[a,b]}$.
In particular, we prove that
$\catC$ is a monoidal categorification of the cluster algebra $K(\catC)$
for $\catC=\catC_\g^0$ and $\catC=\catC_\g^-$.
Note that we need in step (i) above
the existence of
the cluster algebra structure on $K(\catC_\g^{-})$ proved in~\cite{HL16}.

\smallskip

This paper is an announcement whose details will appear elsewhere.

\section{Subcategories} \label{sec: subcategory}
We take the algebraic closure  $\cor$ of $\C(q)$ inside $\bigcup_{m >0}\C((q^{1/m}))$ as the base field for $\uqpg$.
Recall that $\catC_\g$ is the category of finite-dimensional integrable modules over $\uqpg$. There is a family $\set{V(\varpi_i)_c}{i\in I_0,c\in \cor^\times}$ in $\catC_\g$ of simple modules, called the \emph{fundamental modules}.

\smallskip

For simple modules $M$ and $N$ in $\catC_\g$, we say that $M$ and $N$ {\em strongly commute} if $M\tens N$ is simple, and $M$ is \emph{real} if $M^{\tens k}$ is simple for all $k \in \Z_{\ge1}$.

Let us denote by $\Uppsi$ the quiver whose set of vertices is
$$\hI\seteq \{ (i,k) \in I_0 \times \Z \mid k \equiv d(1,i) \ {\rm mod } \; 2 \}.$$
 and
the arrows of $\Uppsi$ consist of two types:
\eq&&\hs{1ex}\left\{
\parbox{50ex}{
\begin{enumerate}[{\rm (A)}]
\item \label{it: arrow1}\hs{2ex}$(i,t) \to (j,s)$ with $d(i,j)=1$ and $s-t=1$,
\item \label{it: arrow2} \hs{2ex}$(i,s+2)\to (i,s)$.
\end{enumerate}
}\right.
\eneq
Here $d(i,j)$ denotes the distance between the vertices $i$ and $j$ in the Dynkin diagram of $\g_0$ and $1\in I_0$ is an arbitrary chosen element.

\smallskip

We say that an infinite sequence
$$\hw_0 = \cdots   s_{i_{-1}}s_{i_0}s_{i_1} \cdots$$ of simple reflections in the Braid group $B(\g_0)$ (\cite{KKOP20A})
of type $\g_0$  is {\em admissible} if
\begin{enumerate} [{\rm (a)}]
\item \label{it:i}
there exists a sequence  $\st{t_k}_{k\in\Z}$ of integers such that
\be[{(1)}]
\item$(i_k,t_k)\in \hI$,
\item $t_{k^+}=t_k+2$, and
\item $t_k>t_{k'}$ if $k>k'$ and $d(i_k,i_{k'})=1$.
\ee
\item \label{it:ii} $s_{i_k}\cdots s_{i_{k+\ell-1}} =w_0$
for all $k \in \Z$,
where $\ell$ denotes the length of longest element $w_0 \in \sfW$.
\end{enumerate}

Here,  for $k \in \Z$ and $j \in I_0$, we set
\begin{equation*}
\begin{aligned}
k^+&\seteq\min\{p \mid k<p ,\; i_k=i_p \}, \\
k^-&\seteq\max\{p \mid    p<k, \; i_k=i_p \},\\
k(j)^+&\seteq\min\{p \mid k \le  p,\; i_p=j \}, \\
k(j)^-&\seteq\max\{p \mid p \le k,\; i_p=j \}.
\end{aligned}
\end{equation*}

\begin{remark}
\bnum
\item We have $i_{k+\ell} =i_k^*$, where
$^*$ denotes the involution on $I_0$ induced by $w_0$.
\item $\hw_0$ completely determines $\st{(i_k,t_k)}_{k\in\Z}$ up to an even translation.
\item For every $k \in \Z$, the reduced expression $s_{i_k}\cdots s_{i_{k+\ell-1}}$ in~\eqref{it:ii}
is \emph{adapted} to some Dynkin quiver $Q$ of type $\g_0$.
Conversely, for any
Dynkin quiver $Q$ of type $\g_0$, there exists a
sequence $\hw_0$ satisfying \eqref{it:i} and \eqref{it:ii}
such that
$s_{i_1}\cdots s_{i_{\ell}}$ is adapted to $Q$.
\ee
\end{remark}
\smallskip

For each $k \in \Z$, we define the fundamental module
$$V[k]^{\hw_0} \seteq V(\varpi_{i_k})_{(-q)^{t_k}}.$$
Then we have
\begin{align*}
V[k^\pm]^{\hw_0}\simeq V[k]^{\hw_0}_{(-q)^{\pm 2}}\qtq
 V[k+\ell]^{\hw_0} = \dual\bl V[k]^{\hw_0}\br,
\end{align*}
where $\dual$ is the right dual functor.
\begin{definition}
For each interval $[a,b]$, we denote by $\catC_\g^{[a,b]}$ the smallest full subcategory of $\catC_\g$  satisfying the following conditions$\colon$
\begin{enumerate}[{\rm (i)}]
\item it is stable under taking subquotients, extensions, tensor products and
\item it contains $V[k]^{\hw_0}$ for all $a \le k \le b$ and the trivial module $\mathbf{1}$.
\end{enumerate}
\end{definition}

\smallskip

\begin{remark} \label{rmk: knownC}
 Many of known subcategories $\shc$ of $\catC_\g$ can be identified with
$\catC_\g^{[a,b]}$ by taking suitable $[a,b]$:
\begin{enumerate}[{\rm (1)}]
\item  $\catC_\g^{[-\infty,\infty]}$ coincides with the subcategory $\catC_\g^0$.
\item  The subcategory $\catC_\cQ$ associated to a
$\rmQ$-data $\cQ$ coincides with
$\catC_\g^{[a,b]}$ for some interval $[a,b]$ with $\cd{[a,b]}\seteq b-a+1=\ell$.
\item By taking $s_{i_1}\cdots s_{i_{\ell}}$ in~\eqref{it:ii}
as adapted to the Dynkin quiver $Q$ with
$ \{1,2\}  \ni  \phi_Q(k) \equiv d(1,i_k)$
(${\rm mod} \; 2$) and   $t_k=\phi(i_k)$ for $1 \le k \le \cd{I_0}$,
 $\catC_N$ can be identified with $\catC_\g^{[a,0]}$  where
$a=1-  (N \times \cd{I_0})$, and
$\catC_\g^-$ can be identified with $\catC_\g^{[-\infty,0]}$.
Those subcategories $\catC_N$, $\catC_\g^-$ of $\catC_\g^0$ are
introduced in \cite{HL10,HL16}.
\end{enumerate}
\end{remark}

\smallskip

\section{Real simple commuting family associated to an admissible chain of $i$-boxes}
Let us fix an admissible sequence $\hw_0$ and $\st{t_k}_{k\in\Z}$.
We write $V[k]$ for $V[k]^{\hw_0}$.
We say that an interval $c=[a,b]$ is an \emph{$i$-box}
if ${i_a}={i_b}$. For each $i$-box $[a,b]$, the module
$M[a,b]$ in $\catC_\g^0$ is defined as follows:
$$M[a,b]\seteq\hd\bl V[b] \tens V[b^-] \tens \cdots \tens V[a^+] \tens V[a] \br,$$
where $\hd(M)$ for $M \in \catC_\g$ denotes the head of $M$. In particular, $V[a]=M[a,a]$.

\smallskip

\begin{theorem} \label{thm: M[a,b]}
\begin{enumerate} [{\rm (i)}]
\item \label{it:A1} $M[a,b]$ is a Kirillov-Reshetikhin module
with a dominant extremal weight $s\varpi_{i_a}$ where $s = \cd{\{ k \mid a \le k \le b, \ i_k=i_a\}}$.
\item \label{it:A2} For $i$-boxes $[a,b]$ and $[c,d]$,
 $M[a,b]$ and $M[c,d]$ commutes if either
$$a^- < c \le d < b^+ \ \ \text{ or } \ \  c^- < a \le b < d^+.$$
\item \label{it:A3}
For any $i$-box $[a,b]$, there exists an exact sequence in terms of $M[a,b]$'s as follows$\colon$
\begin{equation}\label{eq:T-system}
0   \To \hs{-1.7ex} \sotimes_{d(i_a,j)=1}\hs{-.5ex}  M[a(j)^+ ,b(j)^-]  \to    M[a^+ ,b] \tens M[ a,b^-]
 \to M[a,b]\tens M[a^+ ,b^-] \to 0,
\end{equation}
We call it a {\em $T$-system}.
\end{enumerate}
\end{theorem}

\begin{remark} \label{rmk: determinantial module}
For any reduced expression $\rxw=s_{j_1}\cdots s_{j_\ell}$ of $w_0$ (\emph{not necessarily adapted}) and  $[a,b]$ with $j_a=j_b$ and $b-a+1 \le \ell$,
there exists a real simple module $\DC[a,b]^{\rxw}$ in $\shc_{QH}$ of type $\g_0$, called the \emph{determinantial module}, and there exists an exact sequence
(called the T-system)
\begin{align*}
0  \To  \sotimes_{d(i_a,j)=1} \DC[a(j)^+,b(j)^-]\to  \DC[a^+ ,b]   \tens  
\DC[a,b^-]
\to  \DC[a,b]   \tens  \DC[a^+ ,b^-]   \To   0
\end{align*}
in $\shc_{QH}$, which is analogous to~\eqref{eq:T-system}.
More precisely, when $\rxw$ is adapted to some Dynkin quiver $Q$ of type $\g_0$,
\WS functor  $\F_\cQ$ associated with some $\rmQ$-data $\cQ=(Q,\phi_Q)$
transforms the above exact sequence in $\shc_{QH}$
to the T-system \eqref{eq:T-system}.
Thus $M[a,b]$ can be understood as a \emph{quantum affine analogue of the determinantial module}.
(See \cite[Proposition 4.1]{HL11} and \cite{KKKO18} for more detail.)
\end{remark}

\smallskip
For any interval $c\seteq[a,b]$, we introduce $i$-boxes
\begin{align*}
[a,b\}&\seteq[a, b(i_a)^-],  \quad 
\{a,b] \seteq[a(i_b)^+,b]\qtq\\
\LL(c)& \seteq [a-1,b\}, \quad \RR(c) \seteq \{a,b+1].
\end{align*}
\begin{definition} \label{def: admissible chain}
A chain $\frakC$ of $i$-boxes
$$\left( c_k = [a_k,b_k]  \right)_{1 \le k \le l} \quad (l \in \Z_{\ge 1}
 \sqcup \{ \infty \})$$
is called \emph{admissible} if the interval
$\tc_k =[\ta_k,\tb_k] \seteq \bigcup_{1 \le j \le k} [a_j,b_j]$
satisfies $\cd{\tc_k} = k$ and one of the following two statements.
\begin{enumerate}[{\rm (1)}]
\item $c_{k}=\LL(\tc_{k-1})$,
\item $c_k=\RR(\tc_{k-1})$.
\end{enumerate}
\ro Please do not confuse $l$ and $\ell$.\rf\
The sequence of intervals $\st{\tc_k}_{1\le k\le l}$  is called
the {\em envelope} of the chain
$\frakC$ and $\tc_l$ is called the {\em range} of $\frakC$.
\end{definition}

\smallskip

Thus, for an admissible chain $\frakC$ of $i$-boxes, we can associate a pair $(a, \mathfrak{T})$ consisting of an integer $a$
and a sequence $\mathfrak{T} = ( T_1,T_2,\ldots, T_{l-1})$ such that {\rm (i)} $T_i \in \{ \LL,\RR\}$ $(1 \le i \le l-1)$, {\rm (ii)} $a=a_1=b_1$, \\
{\rm (iii)} $\ [\ta_k,\tb_k] = \begin{cases} [\ta_{k-1}-1,\tb_{k-1}] & \text{ if } T_{k-1}=\LL, \\
[\ta_{k-1},\tb_{k-1}+1]  & \text{ if } T_{k-1}=\RR.\end{cases}$\\
Hence we have
$c_k=[a_k,b_k] = T_{k-1}[\ta_{k-1},\tb_{k-1}]$  $(k \ge 2)$,
and the interval $\tc_{k}$ is obtained from $\tc_{k-1}$ by adding an element from the left
or from the right according that $T_{k-1}=\LL$ or $T_{k-1}=\RR$.

\smallskip

For an admissible chain $\frakC=(c_k)_{1\le k\le l}$ with the associated pair $(a, \mathfrak{T})$ and for $1 \le s <l$,
we say that an $i$-box $c_{s}$ is \emph{movable}
if $s=1$ or $T_{s-1} \ne T_s$ $(s\ge2)$.
For a movable $c_s$ in $\frakC$, we define a
new  admissible chain $B_s(\frakC)$ whose associated pair
$(a',\mathfrak{T}')$ is given
\begin{align*}
& {\rm (i)} \ \begin{cases}
a' = a \pm  1 &  \text{if $s=1$ and $T_1 = \RR$ (resp.\ $\LL$),}  \\
a'=a & \text{if $s>1$,}
\end{cases} \\
& {\rm (ii)} \ T'_{k} = T_k \text{ for } k \not\in \{ s-1,s\}, \text{ and } \\
& {\rm (iii)} \ T'_{k}\ne T_k \text{ for } k \in \{ s-1,s\}.
\end{align*}
That is, $B_s(\frakC)$ is the admissible chain obtained from $\frakC$
by moving  $\tc_s$ by $1$ to the right or to the left
inside  $\tc_{s+1}$.

\smallskip

\begin{theorem} \label{thm: real commuting family}
Let $\frakC=(c_k)_{1\le k\le l}$ be an admissible chain and set
$$M(\frakC)\seteq \{ M[a_k,b_k] \mid  1 \le k \le l \}.$$
Then we have the followings:
\begin{enumerate}[{\rm (a)}]
\item \label{it:b} $M(\frakC)$ is a commuting family of real simple modules.
\item \label{it:c} If $M[c,d]$ commutes with all $M[a_k,b_k]$ and $[c,d] \subset [\ta_l,\tb_l]$, then $[c,d] \in \frakC$.
\item \label{it:d} For another admissible chain $\frakC'=(c'_k)_{1\le k\le l}$ with the same range,
there exists a finite sequence
$(t_1,t_2,\ldots,t_r) \hspace{-.3ex} \in \hspace{-.3ex}\{ 1,2,\ldots, l \}^r$ such that
$$   B_{t_r}(\cdots(B_{t_2}(B_{t_1}(\frakC))\cdots)  = \frakC'.$$
\end{enumerate}
\end{theorem}

\section{Monoidal categorification} Let $\K=\Kex \sqcup \Kfr$ be a countable index set. Let $\catC$ be a full subcategory of
$\catC_\g^0$ stable under taking subquotients, extensions and tensor products.

\smallskip

A \emph{monoidal seed in $\catC$} is
a pair $\seed = (\{ M_i\}_{i\in \K },\widetilde B)$ consisting of a commuting family $\{ M_i\}_{i\in\K}$ of real simple objects in
$\catC$ and a $\Z$-valued $\K\times\Kex$-matrix $\widetilde B = (b_{ij})_{(i,j)\in\K\times\Kex}$
such that 
\bnum
\item for each $j \in \Kex$, there exist finitely many $i \in \K$ such that $b_{ij} \ne 0$, 
\item the {\it principal part} $B \seteq (b_{ij})_{i,j \in \Kex}$ is skew-symmetric.
For $i\in\K$, we call $M_i$ the $i$-th {\em cluster variable module} of $\seed$.
\ee
\smallskip

For a monoidal seed  $\seed=(\{M_i\}_{i\in\K}, \widetilde B)$, let $\La^\seed=(\La^\seed_{ij})_{i,j\in\K}$ be the skew-symmetric matrix given by $\La^\seed_{ij}=\Lambda(M_i,M_j)$
(see \cite{KKOP20}).

A monoidal seed $\seed=(\{M_i\}_{i\in\K}, \widetilde B)$ is called
\emph{$\Uplambda$-admissible} if \\
{\rm (i)} $(\Lambda^\seed\tB)_{jk}=-2\delta_{jk}$ for $(j,k)\in\K\times\Kex$, and\\
{\rm (ii)} for each $k\in\Kex$, there exist a simple object $M'_k$ of $\shc$
commuting with $M_i$  for any $i \neq k$
and an exact sequence in $\catC$
\eq\label{eq:mut}&&
0\To\sotimes_{b_{ik} >0}  M_i^{\tens  b_{ik}}\to
M_k \otimes  M_k'\to \sotimes_{b_{ik} <0} M_i^{\tens  (-b_{ik})}\To0.
\eneq

\smallskip

 Under the following two assumptions\\[.5ex]
{\rm (a)} there exists a $\Uplambda$-admissible monoidal seed $\seed=(\{M_i\}_{i\in K},\widetilde B)$ in $\catC$,\\
{\rm (b)} $K(\catC)$ is isomorphic to the cluster algebra $\mathscr{A}([\seed])$,\\[1ex]
KKOP (\cite[Theorem 6.10]{KKOP19C}) proved that
$\catC$ provides a monoidal categorification of
$\mathscr{A}([\seed])$.
Here $[\seed]\seteq\bl[\{M_i]\}_{i\in K},\widetilde B\br$ is a seed in $K(\catC)$,
and $\mathscr{A}([\seed])$ denotes the cluster algebra with the initial seed $[\seed]$.

\smallskip

Set ${\hI}^- \seteq \hI \cap (I_0 \times \Z_{\le 0})$
and let $\Uppsi^-$ be the full subquiver of $\Uppsi$
whose set of vertices is $\hI^-$.
In~\cite{HL16}, Hernandez-Leclerc proved that $\scrA^- \seteq K(\catC_\g^-)$ has a cluster algebra structure whose initial cluster variable modules
$ \{ \M_{(i,t)}  \}_{(i,t) \in \hI^-}$ consist of certain KR-modules.
For a suitable choice of $\hw_0$ (Remark~\ref{rmk: knownC}), we have $\catC_\g^-=\catC_\g^{[-\infty,0]}$ and
$ \{ \M_{(i,t)}  \}_{(i,t) \in \hI^-}$  can be described as $M(\frakC^{-})$  for the following admissible chain $\frakC^{-}$ of
$i$-boxes:
$$ \frakC^- = (0, \frakT=(\LL,\LL,\LL,\ldots) ).$$
More precisely, for $(i,t)=(i_a,t_a)$ ($a\le0$), we have
$$ \M_{(i,t)} =  M[a, 0 \}.$$

The following theorem gives an affirmative answer for the conjecture on $\catC_\g^-$:

\smallskip

\begin{theorem} \label{thm: cg-}
The monoidal seed $$\seed^-\seteq ( M(\frakC^{-}),\tB^-) \text{ is $\Uplambda$-admissible},$$ where $\tB^-$ is the matrix associated with $\Uppsi^-$.
Hence $\catC_\g^-$ provides a monoidal categorification of $\scrA^-$.
\end{theorem}

\smallskip

Now we shall generalize the above theorem to an arbitrary $\catC_\g^{[a,b]}$.

\smallskip

\begin{prop} \label{prop: must be}
Let $\frakC=(c_k)_{1 \le k \le l}$ be an admissible chain of $i$-boxes
with the range $[a,b]$ and the envelope  $\st{\tc_k}_{1\le k\le l}$.
Assume that $\catC\seteq\catC_\g^{[a,b]}$ provides a monoidal categorification of $K(\catC)$ with a $\Uplambda$-admissible monoidal seed
$( M(\frakC),\tB)$.
Let $c_s$ be a movable $i$-box of $\frakC$
and set $\frakC'=B_s(\frakC)$.
\bnum \item
If $\tc_{s+1}\neq c_{s+1}$, then $M(\frakC')$ is equal to $M(\frakC)$
up to a permutation.
\item
If $\tc_{s+1}=c_{s+1}$, then $M(\frakC')$ is the monoidal mutation of
$M(\frakC)$ at $s$. Moreover the corresponding exact sequence
\eqref{eq:mut} is given by the T-system \eqref{eq:T-system}.
\ee
\end{prop}

The above proposition and Theorem~\ref{thm: real commuting family} show that
all $M(\frakC)$ are mutation equivalent.

\smallskip

Now we state our main theorem:

\smallskip

\begin{theorem} \label{thm:main}
For any admissible chain $\frakC=(c_k)_{1\le k \le l}$ for $l \in \Z_{\ge 1} \sqcup \{ \infty \}$
with the range $$\tc_l = [a,b] \quad \text{ for } \ a \le b \in \Z \sqcup \{ \pm\infty \},$$ there exists a $\Uplambda$-admissible monoidal seed $\seed$
of $\catC^{[a,b]}_\g$ such that
\begin{enumerate}[{\rm (i)}]
\item its set of cluster variable modules is $M(\frakC)$,
\item its set of frozen variable modules is
$\{ M[a(i)^+,b(i)^-] \mid i \in I_0, -\infty<a(i)^+\le b(i)^-<+\infty\}$,
\end{enumerate}
and
\begin{enumerate}
\item[{\rm (iii)}] 
$ K(\catC_\g^{[a,b]})$ has the cluster algebra structure with
the initial seed $[\seed]$,
and
$\catC_\g^{[a,b]}$ provides a monoidal categorification of 
$\mathscr{A}([\seed])\simeq K(\catC_\g^{[a,b]})$.
\end{enumerate}
\end{theorem}

By Remark~\ref{rmk: knownC}, we have the following corollary:

\begin{corollary}
The Grothendieck ring $K(\catC_\g^{0})$ has a cluster algebra strucure, and
$\catC_\g^{0}$ provides a monoidal categorification of $K(\catC_\g^{0})$. 
\end{corollary}

\smallskip
\begin{remark} \label{rmk: ext}
We can generalize the above results to an arbitrary quantum affine algebra $U_q'(\g)$
by applying a similar framework with
the results in \cite{KKKO16D,KO18,OT19,OhSuh19}.
\end{remark}

\end{document}